\input amstex
\documentstyle{amsppt}
\input epsf.tex
\input label.def

\def\C{{\Bbb C}}
\def\R{{\Bbb R}}
\def\Z{{\Bbb Z}}
\def\Q{{\Bbb Q}}
\def\N{{\Bbb N}}
\def\Rp#1{\R\roman P^{#1}}
\def\Cp#1{\C \roman P^{#1}}

\def\Cl{\mathop{\roman{Cl}}\nolimits}
\def\la{\langle}
\def\ra{\rangle}

\def\oo{\varnothing}

\def\D{\Delta}
\def\til{\widetilde}
\def\s{\sigma}
\def\e{\varepsilon}
\def\a{\alpha}
\def\b{\beta}
\def\g{\gamma}
\def\d{\delta}

\def\Tors{\roman{Tors}}

\def\sign{\roman{sign}}
\def\SO{\roman{SO}}
\def\U{\roman{U}}
\def\SU{\roman{SU}}

\def\Hom{\operatorname{Hom}}
\def\Sym{\operatorname{Sym}}
\def\l{\lambda}
\def\GC#1#2{G_{#1}(\C^{#2})}
\def\GR#1#2{G_{#1}(\R^{#2})}
\def\tGR#1#2{\til{G}_{#1}(\R^{#2})}
\def\CaR{C_{\a,\R}}
\def\tCaR{\widetilde C_{\a,\R}}
\def\CbR{C_{\b,\R}}

\def\saR{\s_{\a,\R}}
\def\tsaR{\widetilde{\s}_{\a,\R}}
\def\ssaR{s_{\a,\R}}

\def\ovex{={\sssize{}^\pm}\,}
\def\NCp#1{\Cal N^\C_{\boxplus^{#1}}}
\def\NRp#1{\Cal N^\R_{\boxplus^{#1}}}
\def\Nminp#1{\Cal N^{\min}_{\boxplus^{#1}}}
\def\Nep#1{\Cal N^{e}_{\boxplus^{#1}}}

\def\o{\widetilde{\vartheta}}
\def\on{\vartheta}

\def\pii{\pi}

\let\t=\trans

\let\ge\geqslant
\let\le\leqslant

\NoBlackBoxes

 \topmatter
\title
Abundance of $3$-planes on real projective hypersurfaces
\endtitle
\author S.~Finashin, V.~Kharlamov
\endauthor
\address Middle East Technical University,
Department of Mathematics
\endgraf Universiteler Mahallesi,
Dumlupinar Bulvari 1, Ankara 06800 Turkey
\endaddress
\address
Universit\'{e} de Strasbourg et IRMA (CNRS)\endgraf 7 rue Ren\'{e}-Descartes 67084 Strasbourg Cedex, France
\endaddress
\abstract
We show that a generic real projective $n$-dimensional hypersurface
of odd degree $d$, such that
$4(n-2)=\binom{d+3}3$, contains "many" real $3$-planes, namely, in the logarithmic scale
their number has the same rate of growth, $d^3\log d$,
 as the number of complex $3$-planes. This estimate is based on the interpretation of a
suitable signed count of the $3$-planes as the Euler number of an appropriate bundle.
\endabstract
\endtopmatter

\hskip1.8in Everything you can imagine is real.
\vskip3mm
\hskip2.5in  Pablo Picasso
 \vskip6mm


\document
\section{Introduction}

\subsection{Phenomenon of abundance}
Up to our knowledge, the phenomenon of abundance of real solutions in certain real enumerative problems
was observed for the first time in \cite{iks}, where it was shown that the number of real rational curves of degree $d$
interpolating a generic collection
of $3d-1$ real points in the real  projective plane grows, in the logarithmic scale, as fast as the number of complex curves.
Since then, similar abundance phenomena were observed in various other real enumerative problems ({\it cf.}, \cite{iks2}, \cite{GZ}, \cite{KR}).
In particular, in our previous paper \cite{FK} we proved that a generic real
hypersurface of degree $2n-1$
in a real projective space of dimension $n+1$ contains at least
$(2n-1)!!$ real lines, which is approximately  the square root of the number of complex lines (the same bound was obtained by C.~Okonek and A.~Teleman \cite{OT}).
This estimate was based on
the signed counting of the real lines by means  of
the Euler number of an appropriate vector bundle, and as a result gave the following relations
 $$
 n_d^\C\ge n_d^\R\ge n_d^{\R,\min}\ge n_d^e, \quad \log n_d^e\sim \frac12 \log n_d^\C,
 $$
 where $n_d^\C$ and $n_d^\R$ denote respectively the numbers of complex and real lines on a generic real hypersurface
 of odd degree $d$ in a projective space of dimension $\frac{d+3}2$,
 the symbol $n_d^{\R,\min}$ stands for
 the minimum of $n_d^\R$ taken over all such generic hypersurfaces and
$n_d^e$ stands for
the above mentioned signed count of real lines.
(Here, the number $n_d^\R$ depends on the choice of such a hypersurface, while
$n_d^e$, $n_d^{\R,\min}$, and $n_d^\C$ depend only on $d$. A full asymptotic
expansion  of $n_d^\C$ is available due to Don Zagier and can be found in his appendix to \cite{GM}.)

The aim of the present paper is to show that a similar abundance phenomenon holds also as soon
as we count the real $3$-planes on generic real
hypersurfaces of {\it odd} degree $d$ in a real projective space of
an appropriate dimension, which we will still denote by  $n+1$.

To achieve this goal we follow the same approach as in \cite{FK}.
Namely, the variety of (complex or real) $3$-planes on a hypersurface that is
defined by a homogeneous polynomial $f$ of degree $d$ in $n+2$ variables
is viewed as the zero locus $\{s_f=0\}$ in the Grassmannian ($\GC4{n+2}$ or $\GR4{n+2}$, respectively),
where $s_f$ is the determined by $f$
section of the symmetric power $\Sym^d\tau_{4,n+2}^*$ of the dual to the tautological
(complex or real) 4-dimensional vector bundle $\tau_{4,n+2}$ on the corresponding Grassmannian.
It is well-known (see, for example, \cite{DM, Theorem 1.2}) that the section $s_f$ is transversal
to the zero section for a generic choice of $f$.
Thus, if $\dim\{s_f=0\}=4(n-2)-\binom{d+3}3$ is zero, then, in the complex setting,
the Chern number
$c_{4(n-2)}(\Sym^d\tau_{4,n+2}^*)[\GC4{n+2}]$
is equal to the number of complex $3$-planes, while in the real setting,
the Euler number of $\Sym^d\tau_{4,n+2}^*$ on $\GR4{n+2}$ counts the real $3$-planes with signs.

The main feature of such a signed count is its invariance: dependence only on $d$ and not on the choice of a hypersurface (in particular, independence of the topology of the hypersurface). This count is well defined, if $d$ is odd (see Proposition \ref{euler-criterion}), which is
the only interesting case for invariant counting and lower bounds, since in the case of
even $d$
the real locus of the hypersurface
can be empty.
Note also that the count of $3$-planes makes sense only if $n=2+\frac14\binom{d+3}3$ (in higher dimension
the $3$-planes come in families, and in lower dimension their number is zero).

In order to state the results, let us introduce the following notation:
assuming that $d$ is odd,
let us denote by $\Cal N_d^\C$, $\Cal N_d^\R$ and $\Cal N_d^{\R,\min}$ the number of complex $3$-planes on a generic
hypersurface of degree $d$ in a complex projective space of dimension $3+\frac14\binom{d+3}3$, the number of real $3$-planes
on a generic real hypersurface of degree $d$ in a real projective space of dimension $3+\frac14\binom{d+3}3$,
and the minimum of  $\Cal N_d^\R$ taken over all generic
real hypersurfaces as above (the number $\Cal N_d^\C$
does not depend on the choice of a generic hypersurface).
To avoid cumbersome discussions of explicit orientation conventions
needed to fix the sign of the Euler number in question, we take into account only
its absolute value and denote the latter by $\Cal  N_d^e$.
Note that the numbers introduced are linked by the following trivial relations
$$
 \Cal N_d^\C\ge \Cal N_d^\R\ge \Cal N_d^{\R,\min}\ge \Cal  N_d^e\ge 0. \eqtag\label{roughbound}
$$

\theorem\label{main-theorem}
The invariants $\Cal N_d^\C,  \Cal N_d^e$
are positive for each odd $d$ and satisfy the following asymptotic relations as {\rm (}odd{\rm)}
$d\to \infty$:
$$
\log  \Cal N_d^e = \frac1{12}d^3\log d+O(d^3)\le
\log N_d^\C\le \frac16d^3\log{d}+O(d^3) .
$$
\endtheorem

Our conjecture is that, in fact,
$\log \Cal N_d^e\sim \frac12\log \Cal N_d^\C$
which would imply that $\log \Cal N_d^\C\sim\frac16d^3\log d.$
It seems to us that
even the positivity of 
$\Cal N_d^\R$ 
(which follows from $\Cal N_d^e\ne 0$)
was not acknowledged in the literature before.

Amazingly, the answers that we obtain in the real
setting
look more simple than those in the complex setting. Similar phenomena
are observed in other enumerative problems, see \cite{FK} and Section \ref{CR}.

\corollary
 As odd degree $d$ increases, the invariants
$\Cal N_d^\C$ and $\Cal N_d^{\R,\min}$
have the same rate of growth in the logarithmic scale.
More precisely,
$$
\frac1{12}d^3\log d + O(d^3)\le \log \Cal N_d^{\R,\min}\le \log \Cal N_d^\C\le  \frac16d^3\log{d}+O(d^3).
$$
\endcorollary
(For us the same rate of growth
for two sequences $f(n), g(n)$, means existence of a constant $C$
such that $\vert f(n)\vert\le C\vert g(n)\vert$ and $\vert g(n)\vert\le C\vert f(n)\vert$ for all sufficiently big $n$.)

\subsection{Examples, applications, and related results}
In the case of even $d$, we still have $\Cal N_d^\C>0$ (see, for example, \cite{DM, Theorem 2.1}) as well as
$\log N_d^\C\le \frac16d^3\log{d}+O(d^3)$ (see Proposition \ref{complex-upper-bound}). By contrary,
if $d$ is even, then $\Cal N_d^e$ either vanishes or is defined only modulo $2$,
see the explanation in Remark $(1)$ after Proposition \ref{euler-criterion}.

Note that the positivity
$\Cal N_d^\R\ge \Cal N^e_d > 0$ implies that a generic real projective $n$-dimensional hypersurface of odd degree $d$ with
$4(n-2)>\binom{d+3}3$ contains an infinite number of real $3$-planes, and then due to \cite{DM, Theorem 2.1}
 the variety of these real $3$-planes is of (pure) dimension $4(n-2)-\binom{d+3}3$ if $n>6$.
(In fact, O.~Debarre and L.~Manivel  have proved in \cite{DM2}
such positivity and pure dimension results for the variety of real $r$-planes on
odd degree real complete intersections, but only under the assumption that the dimension of the
ambient projective space is large enough. Apparently, for hypersurfaces, due to this dimension assumption their result applies only to $d\le 3$.)

The approach that we develop in this paper can be applied to counting real $(2r-1)$-planes
on real projective $n$-dimensional
 hypersurfaces of any odd degree $d$, under an appropriate dimension condition, that is
$2r(n-2r+2)=\binom{d+2r-1}{2r-1}$. (Counting of even dimensional planes gives a trivial result, since
the dimension of the corresponding vector bundle is odd, and the Euler count, whenever it
gives an integer rather than a modulo 2 residue, would give zero.)
The result of such counting still gives an
invariant, and hence provides, like in the cases $r=1$ (considered in \cite{FK}) and $r=2$ (considered in this paper), an effective universal lower bound for the number of real $(2r-1)$-planes on a hypersurface.
We restricted ourselves here to the case of $3$-planes, since for the moment in the higher dimensional cases we can not
provide an explicit answer, but can only set up an upper bound
(see Subsection \ref{upper_bound}), give an
implicit formula in the form of multivariate
Cauchy integral (see Theorem \ref{general-count}), and suggest a conjecture (see Conjecture \ref{higher-dim-conj}).
\subsection{The content}
In Section 2, we recall
some facts from the complex Schubert Calculus that are related to counting the number of projective subspaces
in  hypersurfaces. In Section 3, we discuss real Schur polynomials and the modifications required to
make similar counting in the real setting. The techniques developed in these sections are applied in Section 4 to prove the main results.
In Section 5 we discuss a few other real enumerative problems that can be solved by using the same methods.

\subsection{Conventions} 
If in a homology or cohomology
 notation the coefficients are not specified, then they are supposed to be integer. The notation $p_i$ for the Pontryagin classes may refer to
 $p_i(\tau_{k,\infty})\in H^*(\GR{k}{\infty})$ as well as for their pull-backs
 in  $H^*(\GR{k}{k+n})$, in $H^*(\tGR{k}{\infty})$, and in $H^*(\tGR{k}{k+n})$
 ($\tilde G$ stands for the Grassmannians of oriented supspaces). This ambiguity should be resolved by the context.

Our decision not to fix explicit orientations results in a number of identities valid up to sign, and we write
$x\ovex{y}$, which means that $x=y$, or $x=-y$.
The symbol \!\!\!$\qed$ marks the end of a proof, or signifies that the corresponding statement is either
a citation or an immediate consequence of previous claims.

\subsection{Acknowledgements} We thank V.~Fock, who helped us to derive the formula
for the Schur coefficients (see Lemma \ref{schur-coefficient}) that is crucial for our calculations.
We thank David Eisenbud for showing how to use the program Macaulay2
to calculate $N_d^\C$ and $N_d^e$ for small values of $d$ (see Example \ref{example-degree3}).
We thank also S.~Fomin, who suggested us some probabilistic heuristic arguments in the favor of the asymptotics
in the problem of counting complex $3$-planes through a given collection of
planes of codimenion $4$
(see Proposition \ref{complex-Catalan}), the asymptotics that we finally justified by applying Laplace's method.

This work was started during visits of the first author to Strasbourg university, continued during visits of the both  authors to the MPIM in Bonn and
to the CIB in Lausanne, and finished during their RIP stay at the MFO in Oberwolfach; we thank all these institutions for providing excellent working conditions.


\section{Elements of Complex Schubert Calculus}

\subsection{Schubert basis}\label{subsection2.1}
By a {\it $k$-partition of $n\in\Z_{\ge0}$} we mean a decreasing integer sequence $\a=(\a_1,\dots,\a_k)$,
$\a_1\ge\dots\ge\a_k\ge0$,
$|\a|=\a_1+\dots+\a_k=n$. Graphically $\a$ is presented as a {\it Young-Ferrers diagram of size $n$}.
For example, {\it the constant $k$-partition} $\bold{m}=(m,\dots,m)$ is presented by the $k\times m$ rectangle.
 In what follows we assume that some $k>0$ is fixed throughout the whole section,
 and omit sometimes the vanishing components of $\a$.

A filtration $0\subset\C^1\subset\C^2\subset\dots$ of $\C^\infty$ yields a CW-decomposition
of Grassmannian $\GC{k}{\infty}$ into open {\it Schubert cells}
$C_\a$ indexed with $k$-partitions,
namely, a $k$-subspace $L\subset\C^{\infty}$ belongs to
$C_\a$ if and only if $\a_{k+1-s}=\min\{j\,|\,\dim(L\cap\C^{j+s})=s\}$,
for each $1\le s\le k$.

With any k-partition $\a$ we associate a homology and a
 cohomology class of Grassmannians as follows.
The closure $\Cl(C_\a)$ is the so-called {\it Schubert variety};
being equipped with the complex orientation it yields the {\it Schubert class}
$[C_\a]\in H_{2n}(\GC{k}{\infty})$, where $n=|\a|$.
The cohomology class $\s_\a\in H^{2n}(\GC{k}{\infty})$ associated to $\a$ is characterized by
$$
\s_\a([C_\b])=\cases
1 \text{ if $\a=\b$},\\
0 \text{ if $\a\ne\b$.}
\endcases
$$
In other words, the classes $\s_\a$ taken over all $k$-partitions of size $n$ form an additive basis in
$H^{2n}(\GC{k}{\infty})$ such that $h=\sum_\a ({h([C_\a])}\s_\a$ for any $h\in H^*(\GC{k}{\infty})$.

Note that the Schubert cell $C_\a\subset\GC{k}{\infty}$ is contained in a finite dimensional Grassmannian
$\GC{k}{k+m}$ if and only if $\a_1\le m$,
that is if the Young-Ferrers diagram of $\a$ lies inside the $(k\times m)$-rectangle diagram. It follows that the additive bases
of $H_*(\GC{k}{k+m})$ and $H^*(\GC{k}{k+m})$ are given respectively
by $[C_\a]$ and $\s_\a$ such that $\a_1\le m$
(here and below, to simplify notation, we denote by $\s_\alpha$ not only the class in $H^{*}(\GC{k}{\infty})$
but also its pull-back in $H^*(\GC{k}{k+m})$).

We say that $k$-partitions
$\a$ and $\b$ are {\it m-complementary to each other} if
$\a_i+\b_{k+1-i}=m$ for $i=1,\dots,k$ (so that, in particular, $\a_1,\b_1\le m$).
It is well-known (and easy to check) that the Schubert cycles of $m$-complementary $k$-partitions
are Poincare-dual in $\GC{k}{k+m}$.

\proposition\label{complex-duality}
Schubert classes
$[C_\a]$ and $[C_\b]$ have intersection index $1$ in $\GC{k}{k+m}$ if k-partitions $\a$ and $\b$ are
m-complementary, and index $0$ if not.
\qed\endproposition

\subsection{Schur polynomials} Denote by $U_k$ the unitary group and by $U_1^k$ its
maximal torus formed by the diagonal matrices.
The inclusion map $U_1^k\subset U_k$ induces
a map of the classifying spaces
$$\CD
BU_1^k=(\Cp\infty)^k @>\phi>> \GC{k}{\infty}=BU_k,
\endCD$$
and the induced cohomology map $\phi^*\:H^*(\GC{k}{\infty})\to H^*(\Cp\infty)^k\cong\Z[z_1,\dots,z_k]$
is independent of the choice of the maximal torus, since such tori form a single conjugacy class.

For the unitary groups, the so called {\it splitting principle} can be expressed formally as follows.

\theorem\label{sym-ring} {\rm [\, E.g., \cite{BT}\,]}
The ring homomorphism $\phi^*$ is monomorphic and its image coincides with the subring $\Z^S[z]\subset \Z[z]$ formed
by the symmetric polynomials. The Chern classes $c_r(\tau_k^*) \in H^*(\GC{k}{\infty})$, $1\le r\le k$,
of the dual tautological vector bundle $\tau_k^*$ over $\GC{k}{\infty}$
are sent to the elementary symmetric polynomials
$$
\phi^*(c_r)=\e_r(z_1,\dots,z_k)=\sum_{i_1<\dots<i_k}z_{i_1}\dots z_{i_k}.
$$
These classes $c_r=c_r(\tau^*_k)$
are multiplicative generators of the ring $H^*(\GC{k}{\infty})$.
\qed\endtheorem

We call $z_i$, $1\le i\le k$, the {\it Chern roots} as they are the formal roots of
$t^k-c_1t^{k-1}+\dots+(-1)^kc_k$. For each $h\in H^*(\GC{k}{\infty})$, we
say that $\phi^*(h)\in\Z^S[z]$ is the {\it root polynomial} representing class $h$.

The root polynomials $s_\a\in\Z^S[z]$ that represent classes $\s_\a\in H^{2n}(\GC{k}{\infty})$, $n=|\a|$,
are called {\it Schur polynomials}.
The relation $h=\sum_\a ({h([C_\a])}\s_\a$ implies
$$
\phi^*(h)=\sum_\a \l_\a(h) s_\a,
$$
where $\l_\a(h)=h[C_\a]\in\Z$ will be called {\it the Schur coefficients} for $h$ (or for $\phi^*(h)$).

Recall that the Schur polynomials can be calculated using the following
{\it generalized Vandermonde polynomial}
$$V_{\a+\d}(z)=\sum_{\t\in S_k}\sign(\t)z_{\t(1)}^{\a_1+k-1}z_{\t(2)}^{\a_2+k-2}\dots z_{\t(k)}^{\a_k}$$
where $\sign(\t)$ is the sign of a permutation $\t\in S_k$,
and $\d=(k-1,k-2,\dots,1,0)$.
 Recall that the usual Vandermonde  polynomial is
 $$V_\d(z)=\sum_{\t\in S_k}\sign(\t)z_{\t(1)}^{k-1}z_{\t(2)}^{k-2}\dots z_{\t(k-1)}=\prod_{1\le i<j\le k}(z_i-z_j).$$

\proposition\label{schur-vandermond} {\rm [\, E.g., \cite{St}, Theorem 7.15.1\,]}
For any $k$-partition $\a$ we have
$$
\hskip2in s_\a=\frac{V_{\a+\d}}{V_\d}.\hskip2in \qed$$
\endproposition

\remark{Example}\label{schur-examples}
\roster\item
The Schur polynomial $s_{1,\dots,1}$ with $r\le k$ components ``$1$'' is the
elementary symmetric polynomial
$$
\e_r(z_1,\dots,z_k)=\sum_{i_1<\dots<i_r}z_{i_1}\dots z_{i_r}=\phi^*(c_r),
$$
it is the root polynomial of $c_r=\s_{1,\dots,1}$ ({\it cf.}, Theorem \ref{sym-ring}).
\item
The Schur polynomial
$s_{m,\dots,m}$ with $k$ components ``$m$'' equals $(z_1\dots z_k)^m$, it is the root polynomial of $c_k^m$.
\endroster
\endremark

\subsection{Multivariate Cauchy integral formula for the coefficients $\l_\a$}
\lemma\label{schur-coefficient}
For any $f\in \Z^S[z]$, $z=(z_1,\dots,z_k)$, and any $k$-partition $\a$, we have
$$
\l_\a(f)=\frac1{k!(2\pi i)^k}\int_{T^k}f(z)\overline{s_\a}(z)V_\d(z)\overline{V_\d}(z)\,\frac{{\roman d}z}{z},
$$
where $T^k=\{z\in\C^k\,:\,|z_1|=\dots=|z_k|=1\}$ and
$\frac{{\roman d}z}{z}=\frac{{\roman d}z_1}{z_1}\dots\frac{{\roman d}z_k}{z_k}$.
\endlemma

\proof
The monomials $z^\a$ with  $\a=(\a_1,\dots,\a_k)$, $\a_i\ge 0$,
form an orthonormal basis in $\C[z]$ with respect to the inner product
$$\la f_1,f_2\ra=\frac1{(2\pi i)^k}\int_{T^k}f_1(z)\overline{f_2(z)}\frac{{\roman d}z}{z}.$$
It follows that for each pair of $k$-partitions $\a$ and $\b$,
$\la V_{\a+\d},V_{\b+\d}\ra=k!\la z^{\a+\d},z^{\b+\d}\ra$.
Thus, according to Proposition \ref{schur-vandermond},
the Schur polynomials $s_\a$ form an orthonormal basis in the
vector space $\Z^S[z]\otimes\C$, with respect to the modified inner product
$$
\la f_1,f_2\ra_{\Sym}=
\frac1{k!(2\pi i)^k}\int_{T^k}f_1(z)\overline{f_2(z)}V_\d(z)\overline{V_\d(z)}\frac{{\roman d}z}{z},\eqtag\label{schur-vander}
$$
namely,
$$\aligned
\la s_\a,s_\b\ra_{\Sym}&=\la\frac{ V_{\a+\d}}{V_\d},\frac{V_{\b+\d}}{V_\d}\ra_{\Sym}=
\frac1{k!(2\pi i)^k}\int_{T^k}V_{\a+\d}(z)\overline{V_{\a+\d}(z)}\frac{{\roman d}z}{z}\\
&=\frac1{k!}\la V_{\a+\d},V_{\b+\d}\ra=\la z^{\a+\d},z^{\b+\d}\ra.
\endaligned$$
The claim of the Lemma follows now from $\l_\a(f)=\la f,s_\a\ra_{\Sym}$ and (\ref{schur-vander}).
\endproof

\corollary\label{fund-class-evaluation}
For each $h\in H^{2km}(G_k(\C^{\infty}))$ its value $h([G_k(\C^{k+m})])$
on the fundamental class
$[G_k(\C^{k+m})]\in H_{2km}(G_k(\C^{\infty}))$ is given by
$$h([G_k(\C^{k+m})])=\l_{m,\dots,m}(\phi^*(h))=
\frac1{k!(2\pi i)^k}\int_{T^k}\frac{\phi^*(h)(z)}{z^{\bold{m}}}V_\d(z)\overline{V_{\d}}(z)\,\frac{{\roman d}z}{z},$$
where $\phi^*(h)\in\Z^S[z]$ is the root polynomial of $h$
and $z^\bold{m}$ in the denominator stands for $(z_1\dots z_k)^m$.
\endcorollary

\proof
We apply Lemma \ref{schur-coefficient} to $\a=(m,\dots,m)$ and use Example \ref{schur-examples}(2).
\endproof

\subsection{Counting $(k-1)$-planes on projective hypersurfaces}\label{(k-1)-planes}
We denote by  $c_{\roman{top}}\in H^{*}(G_k(\C^{k+m}))$ the top Chern class of the
symmetric power $\Sym^d(\tau^*_{k,m})$ of the dual to the tautological bundle $\tau^*_{k,m}$
on $G_k(\C^{k+m})$,
and by $f_d(z)=\phi^*(c_{\roman{top}})\in\Z^S[z]$ the
root polynomial of $c_{\roman{top}}$.
The number of $(k-1)$-planes in a projective hypersurface
can be evaluated as the following Chern number.

\proposition\label{chern-number-interpretation} {\rm [\, E.g., \cite{DM}\,]}
Assume that $X\subset P^{m+k-1}$ is a generic hypersurface of degree $d\ge1$, and that $mk=\binom{d+k-1}{k-1}$.
Then, $X$ contains a finite number of projective $(k-1)$-planes
and this number is equal to
$c_{\roman{top}}$ evaluated on the fundamental class $[G_k(\C^{k+m})]$.
\qed\endproposition
Let us denote the above Chern number $c_{\roman{top}}([G_k(\C^{k+m})])$
by $\Cal N_{d,k}^\C$.
Proposition \ref{chern-number-interpretation} together with Corollary \ref{fund-class-evaluation} provides the following integral formula.

\corollary\label{N_d-via-integral}
If $mk=\binom{d+k-1}{k-1}$, then
$$
\Cal N_{d,k}^\C=
\l_{m,\dots,m}(f_d)=
\frac1{k!(2\pi i)^k}\int_{T^k}\frac{f_d(z)}{z^{\bold{m}}}V_\d(z)\overline{V_{\d}}(z)\,\frac{{\roman d}z}{z}.
\qed
$$
\endcorollary

To estimate this integral, we need the following well-known {\it root factorization formula} for $f_d$.

\proposition\label{f_d-formula} {\rm [\, E.g., \cite{DM}\,]}
For any $d\ge1$,
$$f_d(z)=\prod_{\smallmatrix
\ell_1+\dots+\ell_k=d\\ \ell_1,\dots, \ell_k\ge 0\endsmallmatrix}
(\ell_1z_1+\dots+\ell_kz_k).\qed
\eqtag\label{root-factorization}$$
\endproposition
We call the factors in the
right hand side of (\ref{root-factorization}) the {\it root factors}.

\subsection{An upper bound}\label{upper_bound}
In this section we study the growth rate of the sequence $\Cal N_{d,k}^\C$ in the logarithmic scale.

\proposition\label{complex-upper-bound}
If we fix $k\ge1$ and vary $d\ge1$ so that
$\frac1k\binom{d+k-1}{k-1}$ is integer, then
the invariant $\Cal N_{d,k}^\C$ defined in Subsection \ref{(k-1)-planes}
has the following 
asymptotic upper bound
as $d\to\infty${\rom :}
$$
\log(\Cal N_{d,k}^\C)\le  \frac1{(k-1)!}d^{k-1}\log d+ O(d^{k-2}\log d).$$
\endproposition

\proof Since $|l_1x_1+\dots+ l_kx_k|\le \ell_1+\dots +\ell_k=d$ at each point of $T^k$, the integral formula given in Corollary \ref{N_d-via-integral}
implies that there exists a constant $C$ such that $\Cal N_{d,k}^\C\le Cd^{b(d)}$,
where $b(d)=\binom{d+k-1}{k-1}=\frac1{(k-1)!}d^{k-1}+ O(d^{k-2}).$
\endproof

\subsection{Conjecture}\label{higher-dim-conj} Some heuristic arguments make plausible to conjecture that the sequence $\Cal N_{d,k}^\C$
 has the following logarithmic asymptotics:
$$\log(\Cal N_{d,k}^\C)\sim \frac1{(k-1)!}d^{k-1}\log d\quad \text{for $k$ fixed and $d\to \infty$}.$$


\section{Elements of Real Schubert Calculus}

\subsection{Orientability and the Euler class for the symmetric powers}
We denote the tangent bundle of a real Grassmannian $G_k(\R^{k+m})$ by
$T_{k,m}$ and the tautological $k$-dimensional vector bundle over $G_k(\R^{k+m})$ by $\tau_{k,m}$.

\lemma\label{Gr-orientability} For any $k,m\ge0$, we have
$$
w_1(T_{k,m})=(k+m)w_1(\tau_{k,m}).
$$
 In particular, $G_k(\R^{k+m})$ is orientable if and only if $k+m$ is even.
\endlemma

\demo{Proof} Note that $T_{k,m}=\Hom(\tau_{k,m},\tau_{k,m}^\perp)$,
where $\tau_{k,m}^\perp$ is the $m$-dimensional vector bundle
orthogonal to $\tau_{k,m}$, so that $\tau_{k,m}+\tau_{k,m}^\perp$
is a trivial bundle and, in particular, $w_1(\tau_{k,m})=w_1(\tau_{k,m}^\perp)$.
 Following the splitting principle, we write
$w_1(\tau_{k,m})=\sum_{i=1}^ka_i$ and $w_1(\tau_{k,m}^\perp)=\sum_{j=1}^{m}b_j$,
where $a_i,b_j\in H^1(G_k(\R^{k+m};\Z/2))$,
which gives
$$
w_1(T_{k,m})=\sum_{i,j}(a_i+b_j)=mw_1(\tau_{k,m})+kw_1(\tau_{k,m}^\perp)=(k+m)w_1(\tau_{k,m}). \quad \qed $$
\enddemo

\lemma\label{Sym-orientability}
\roster\item The vector bundle $\Sym^d(\tau^*_{k,m})$ is orientable if and only if $\binom{d+k-1}{k}$ is even.
 \item
If $\dim\Sym^d(\tau^*_{k,m})=\dim G_k(\R^{k+m})$, then
$w_1(\Sym^d(\tau^*_{k,m}))={dm}\,w_1(\tau_{k,m})$. In particular, under
the assumption that $\dim\Sym^d(\tau^*_{k,m})=\dim G_k(\R^{k+m})$
the bundle $\Sym^d(\tau^*_{k,m})$ is
orientable if and only if $dm$ is even.
\endroster
\endlemma

\proof For proving (1), we use the splitting principle and obtain the following expression for
the total Stiefel-Whitney class, $W_*$, of the symmetric power $\Sym^d(\tau_{k,m})$:
$$W_*(\Sym^d(\tau_{k,m}))=\prod_{\smallmatrix\ell_1+\dots+\ell_k=d\\ \ell_1,\dots, \ell_k\ge 0\endsmallmatrix}(1+\ell_1a_1+\dots+\ell_ka_k)
$$
with respect to $w_1(\tau_{k,m})=\sum_{i=1}^ka_i$. Hence,
$$
w_1(\Sym^d(\tau_{k,m}))= \sum_{\smallmatrix\ell_1+\dots+\ell_k=d\\ \ell_1,\dots, \ell_k\ge 0\endsmallmatrix}(\ell_1a_1+\dots+\ell_ka_k)=n(a_1+\dots+a_k)
$$
where
$$
n=\frac{d}k\binom{d+k-1}{k-1}=\binom{d+k-1}{k},
$$
from where $w_1(\Sym^d(\tau_{k,m}))=0$ if and only if $n=\binom{d+k-1}{k}$ is even.

To deduce (2), note that $\binom{d+k-1}{k-1}=\dim\Sym^d(\tau_{k,m})=\dim G_k(\R^{k+m})=km$ implies
$n=\frac{d}kkm=dm$.
\endproof

As an immediate consequence we obtain the following result.

\proposition\label{euler-criterion}
Assume that the dimension of the vector bundle
$\Sym^d(\tau^*_{k,m})$, that is $\binom{d+k-1}{k-1}$, is equal to the dimension $km$ of the Grassmannian
$G_k(\R^{k+m})$. Assume also that $k+m=dm\mod2$.
Then, the Euler number
$e(\Sym^d(\tau^*_{k,m}))[G_k(\R^{k+m})]$ with respect to the local coefficient system twisted by $w_1(G_k(\R^{k+m}))$
is an integer well-defined up to sign.
\qed\endproposition

\remark{Remarks}
\roster\item
The first assumption of Proposition \ref{euler-criterion} is always satisfied in what follows, since it simply
signifies that the virtual dimension of the variety $F_{k-1}(X)$ of $(k-1)$-planes contained in
a hypersurface $X\subset P^{k+m-1}$ of degree $d$ is $0$.
Note that for even $d$ there exist real hypersurfaces $X$ with $X(\R)=\oo$, for instance,
the Fermat hypersurface, and thus, a signed count of the real $(k-1)$-planes
on such hypersurfaces (if invariant) gives $0$.
Note also that if $km$ is odd, then even when the Euler number is well defined, it is
equal to zero, as it happens for any real odd dimensional vector bundle.

For odd $d$ the second assumption just means that $k$ is even, so, the case of even $k$ and odd $d$ is the
only one which makes sense to study.
 As the case $k=2$ was analyzed in \cite{FK}, we are concerned in what follows
with the case $k=4$.
\item
Calculation of the Euler number $e(\Sym^d(\tau^*_{k,m}))[G_k(\R^{k+m})]$
can (and will) be done by passing to the Grassmannian
$\til{G}_k(\R^{k+m})$ of orientable
$k$-planes and its (orientable)
dual tautological bundle $\til{\tau}^*_{k,m}$, because of the relation
$$e(\Sym^d(\tau^*_{k,m}))[G_k(\R^{k+m})]\ovex\frac12e(\Sym^d(\til{\tau}^*_{k,m}))[\til{G}_k(\R^{k+m})].
\eqtag\label{count_on_oriented_grassmannian}
$$
\endroster\endremark

\subsection{Pontryagin classes and the real root polynomials}
Let $\tGR{2k}{\infty}$ be the Grassmannian of oriented $2k$-planes in $\R^\infty$,
and  $\GR{2k}{\infty}$ be that of non oriented  $2k$-planes.
Denote by ${\o}\:\tGR{2k}{\infty}\to\GC{2k}{\infty}$ the composition of the double covering
$\pi\:\tGR{2k}{\infty}\to \GR{2k}{\infty}$ and of the inclusion
$\on\: \GR{2k}{\infty}\subset  \GC{2k}{\infty}$.
The following description of the integer cohomology ring $H^*(\tGR{2k}{\infty})/\Tors$ is classical.

\theorem\label{EP-ring}{\rm [\,E.g., Brown \cite{Br}\,]}
{\tenrm(1)} The ring $H^*(\GR{2k}{\infty})/\Tors$ is freely generated by the
Pontryagin classes
$p_i=(-1)^i \on^*(c_{2i})\in H^{4i}(\GR{2k}{\infty}), 1\le i\le k .$

{\tenrm(2)} The ring $H^*(\tGR{2k}{\infty})/\Tors$ is generated
by the Pontryagin classes
$$\tilde p_i=\pi^*(p_i)=(-1)^i\o^*(c_{2i})\in H^{4i}(\tGR{2k}{\infty}),\ \  1\le i\le k$$
and the Euler class $e_{2k}\in H^{2k}(\tGR{2k}{\infty})$,
the only defining relation in these generators is
$\tilde p_k=e_{2k}^2$.
\qed\endtheorem

In the rest of the paper, we allow ambiguity and keep traditional notation $p_i$ for the classes $\tilde p_i$, and moreover,
use the same notation for the Pontryagin classes induced in the Grassmannians
$H^*(\GR{2k}{n})/\Tors$ and $H^*(\tGR{2k}{n})/\Tors$.

The embedding $\SO_2^k\subset \SO_{2k}$ given by $2\times 2$ special orthogonal matrices along the diagonal, gives a maximal torus
in $\SO_{2k}$ and induces a map between the classifying spaces
$$\CD
(\Cp{\infty})^k=(B\SO_2)^k @>\phi_\R>> B\SO_{2k}=\tGR{2k}{\infty}.
\endCD$$
This map associates to a k-tuple $(\xi_1,\dots,\xi_k)$ of $\SO_2$-bundles their Whitney sum
$\xi_1\oplus\dots\oplus\xi_k$. We denote by $x_1,\dots, x_k$ the standard (Euler class) generators of $H^*((B\SO_2)^k)$
and, for a given a class $h\in H^*(\tGR{2k}{\infty})$, by a {\it real root polynomial} of $h$
we mean the polynomial $\phi_\R^*(h)\in\Z[x_1,\dots,x_k]=H^*((B\SO_2)^k)$.

\lemma\label{real-root-pullback}
Consider any class $h\in H^*(\GC{2k}{\infty})$ and let $h_\R=\o^*(h)\in H^*(\tGR{2k}{\infty})$ denote its pull-back. Then the real root polynomial $\phi_\R^*(h_\R)(x)$ is obtained from the complex one,
$\phi^*(h)(z)$, by letting $z_{2k-1}=-z_{2k}=x_k$, that is
$$
\phi_\R^*(h_\R)(x_1,\dots,x_k)=\phi^*(h)(x_1,-x_1,\dots,x_k,-x_k).
$$
\endlemma

\proof
The tautological embedding $\SO_2\to\SU_2\subset\U_2$ is conjugate to the embedding
 $\SO_2\to\SU_2\subset \U_2$ defined as the composition of the {\it antidiagonal homomorphism}
$\SO_2\to \SO_2\times\SO_2$, $g\mapsto (g,g^{-1})$, the
isomorphism $\SO_2\times\SO_2\cong\U_1\times\U_1$ and the maximal torus inclusion $\U_1\times\U_1\subset U_2$.
This leads to a commutative up to conjugation diagram
$$
\CD
(\SO_2)^k @>>> (\U_1)^{2k}\\
@VVV  @VVV\\
\SO_{2k} @>>> \U_{2k}
\endCD
$$
with the ``vertical'' block-diagonal inclusion maps, the k-th power of the antidiagonal homomorphism
$\SO_2\to\SO_2\times\SO_2\cong \U_1\times\U_1$ at the ``top'', and the tautological homomorphism at the ``bottom''.
 Passing to the corresponding classifying spaces we obtain a commutative up to homotopy diagram
$$
\CD
(B\SO_2)^k=&(\Cp\infty)^k @>(\D_a)^k>> &(\Cp\infty)^{2k}&=(B\U_1)^{2k}\\
 &@V\phi_\R VV  &@VV\phi_\C V\\
B\SO_{2k}=&\tGR{2k}{\infty} @>\o>> &\GC{2k}{\infty}&=B\U_{2k}
\endCD
$$
where $\D_a\:\Cp\infty\to (\Cp\infty)^2$ is the {\it antidiagonal map} $z\mapsto (z,\bar z)$.
The induced map in cohomology yields a commutative diagram
$$
\CD
H^*(\Cp\infty)^k=&\Z[x_1,\dots,x_k] @<(\D_a^*)^k<< \Z[z_1,\dots,z_{2k}]&=H^*(\Cp\infty)^{2k}\\
&@A\phi_\R^* AA  @AA\phi_\C^* A\\
&H^*(\tGR{2k}{\infty}) @<<\o^*< H^*(\GC{2k}{\infty})
\endCD
$$
where $z_{2i-1}\mapsto x_i$, $z_{2i}\mapsto -x_i$, $i=1,\dots, k$.
\endproof

We consider
two subrings of the ring $\Z[x]$, $x=(x_1,\dots, x_k)$.
The {\it Pontryagin ring}
$\Z^P[x]=\Z[x_1^2,\dots,x_k^2]\cap\Z^S[x]$ is formed by the symmetric polynomials in $x_i^2$, $1\le i\le k$, and
the {\it Euler-Pontryagin ring}
$\Z^{EP}[x]=\Z[x_1^2,\dots,x_k^2,x_1\dots x_k]\cap\Z^S[x]$, where
$\Z[x_1^2,\dots,x_k^2,x_1\dots x_k]$ is generated by the
squares $x_i^2$ and the product $x_1\dots x_k$.

\proposition\label{real-root-map} For any $k\ge1$ and $x=(x_1,\dots,x_k)$:
\roster\item
 the map $\phi_\R^*$ yields a monomorphism
$H^*(\tGR{2k}{\infty})/Tors\to \Z[x]$;
\item
 the image of $\phi_\R^*$ is $\Z^{EP}[x]$;
\item
 the image of $\phi_\R^*\circ\o^*$ is the image of  $\phi_\R^*\circ\pi^*$ and is $\Z^{P}[x]$;
\item
  $\phi_\R^*(p_i)=\e_i(x_1^2,\dots,x_k^2)$ and $\phi_\R^*(e_{2k})=x_1\dots x_k$.
\endroster\endproposition

\proof
Follows immediately from Lemma \ref{real-root-pullback}, Theorem \ref{EP-ring} and Theorem \ref{sym-ring}.
\endproof

\corollary\label{P_and_EP-rings}
The maps  $\phi_\R^*\circ\pi^*$ and $\phi_\R^*$ identify the rings
$H^*(\GR{2k}{\infty})/Tors$ and $H^*(\tGR{2k}{\infty})/Tors$ with $\Z^{P}[x]$ and $\Z^{EP}[x]$,
respectively. The map
$$\pi^*\:H^*(\GR{2k}{\infty})/Tors\to H^*(\tGR{2k}{\infty})/Tors$$
induced by the projection $\pi\:\tGR{2k}{\infty}\to \GR{2k}{\infty}$
 is identified then with the inclusion homomorphism $\Z^{P}[x]\to \Z^{EP}[x]$.\quad \qed
\endcorollary

\subsection{Integral Schubert cycles of infinite order in $\GR{2k}{\infty}$}
As in the complex case, 
$k$-partitions $\beta$ determine a CW-decomposition of $\GR{k}{\infty}$ into {\it real Schubert cells} $\CbR=C_\b\cap \GR{k}{\infty}$ with the property
$\CbR\subset\GR{k}{k+m}$ if and only
if $\b_1\le m$. The cells $\CbR$ lead in general only to $\Z/2$-homology classes
$$
[\CbR]_2\in H_{|\b|}(\GR{k}{k+m};\Z/2)\to H_{|\b|}(\GR{k}{\infty};\Z/2), \quad \b_1\le m\le\infty
$$
(here and further on, the same notation is used
for the classes  in Grassmannians for different $m$, including $m=\infty$).

For a $k$-partition $\b=(\b_1,\dots,\b_k)$,
denote
by $2\b$ the $k$-partition $(2\b_1,\dots,2\b_k)$
and by $\b(2)$ the $2k$-partition $(\b_1,\b_1,\dots,\b_k,\b_k)$, where each of the $\b_i$ is repeated twice.
We say that $2k$-partition $\a$ is an {\it even partition} if it can be presented as
$\a=2\b(2)$, and that $\a$ is an {\it odd partition} if it has form
$\a=2\b(2)+\bold{1}=(2\b_1+1,\dots,2\b_k+1)$ for some $k$-partition $\b$.

If $\a$ is an even $2k$-partition, then
the real Schubert cell $\CaR$ yields an integral homology class, which we denote
$$[\CaR]\in H_{|\a|}(\GR{2k}{2k+m})/\Tors,\ \ \a_1\le m\le\infty.
$$
We postpone our choice of the sign of $[\CaR]$ until stating
Proposition \ref{Fuks-property}.

\theorem\label{basis-pontryagin-nonorientable}{\rm [\,Pontryagin \cite{P}\,]}
Assume that $m\ge0$ is even and $r$ or $2km-r$ is less than $m$.
Then the Schubert classes
$[\CaR]$ for all even $2k$-partitions $\a$ of dimension $|\a|=r$ form
a basis in the group $H_r(\GR{2k}{2k+m})/\Tors$.
\qed\endtheorem

\remark{Remark}
Pontryagin's claim (see \cite{P}, Theorem 1) concerns the orientable Grassmannians, which is
stronger than what we stated here. On the other hand, his claim
covers only the case of $2km-r<m$. The case $r<m$ can be derived from that one via Poincare duality,
since the Poincare-dual Schubert cycles are represented by $m$-complementary even $2k$-partitions (for even $m$).
Indeed,
for even $m$ (in which case Grassmannian $\GR{2k}{2k+m}$ is orientable by Lemma \ref{Gr-orientability}),
the intersection index of Schubert cycles $[\CaR]\circ[C_{\b,\R}]$ in $\GR{2k}{2k+m}$
is $\pm1$ if $\a$ and $\b$ are complementary even $2k$-partitions, and otherwise that index is $0$.
\endremark

\subsection{Real Schur polynomials}
Letting $m\to\infty$ in Theorem \ref{basis-pontryagin-nonorientable}, we can conclude that classes
$[\CaR]$ for all even $2k$-partitions $\a$ form an additive basis in $H_*(\GR{2k}{\infty})/\Tors$.
Let us introduce
the dual additive basis, $\{\saR\}_{\text{even } \a}$, in $H^*(\GR{2k}{\infty})/\Tors$;
namely, let $\saR$ be
the integral cohomology class taking value $1$ on $[\CaR]$ and
vanishing on the other
integral Schubert classes.
Then, we take the pull-back of $\saR$ by
the double covering $\pii\:\tGR{2k}{\infty}\to\GR{2k}{\infty}$,
that is
$\tsaR=\pii^*(\saR)\in H^*(\tGR{2k}{\infty})$,
and define a
{\it real Schur polynomial of an even $2k$-partition $\a$}
as the root polynomial
$$
\ssaR=\phi_\R^*(\tsaR)\in\Z^{P}[x]\subset\Z^{EP}[x],
$$
({\it cf.}, Corollary \ref{P_and_EP-rings}).

The relation between the real and the complex Schur polynomials
can be specified by an appropriate choice of orientations as follows.

Consider the ring isomorphism $T$ between the Chern ring
$$
\Z^{S}[z_1,\dots,z_k]\cong H^*(\GC{k}{\infty})
$$ and
the Pontryagin ring
$$\Z^{P}[x_1,\dots,x_k] \cong H^*(\GR{2k}{\infty})/\Tors
$$
that sends, in terms of polynomial rings, $z_i$ to $x_i^2$, 
or
 equivalently, the Chern classes $c_i$ to the Pontryagin classes $p_i$, $i=1,\dots, k$,
 in terms of cohomology rings.

\proposition\label{Fuks-property}{\rm [\,E.g., Fuks \cite{Fuks}, Section 2.3.E\,]}
There exist such a choice of orientations of $C_{\alpha,\R}$ for all even $2k$-partitions $\a=2\b(2)$ that 
$c_1^{\g_1}\dots c_k^{\g_k}[C_\b]=p_1^{\g_1}\dots p_k^{\g_k}[\CaR]$
for any Chern class $c_1^{\g_1}\dots c_k^{\g_k}$. The
isomorphism $T\: \Z^{S}[z_1,\dots,z_k]\to \Z^{P}[x_1,\dots,x_k]$
sends the Schur polynomial
$s_\b$ to the real Schur polynomial $\ssaR$ defined under this choice.
\qed\endproposition

In what follows we suppose that for each even partition $\a$ the orientation of $C_{\alpha,\R}$  is fixed in accordance with
Proposition \ref{Fuks-property}.

 Now we extend the definition of the real root polynomials to odd $2k$-partitions,
 $\a'=\a+\bold{1}$ (where  $\a$ is an even partition)
 by letting $\widetilde{\s}_{\a',\R}$ be the product of
 $\tsaR$ by the Euler class $e_{2k}\in H^{2k}(\tGR{2k}{\infty})$ and define similarly the
{\it real Schur polynomial } of $\a'$ as the root polynomial
$$s_{\a',\R}=\phi_\R^*(\widetilde{\s}_{\a',\R})\in\Z^{EP}[z].
$$

\corollary\label{real-complex-Schur}
For any even  $2k$-partition $\a=2\b(2)$, and the odd partition $\a'=\a+\bold{1}$
we have
\roster\item
$\ssaR(x_1,\dots,x_k)=s_{\b}(x_1^2,\dots,x_k^2);$
\item
$\widetilde{s}_{\a',\R}(x_1,\dots,x_k)=x_1\dots x_ks_{\b}(x_1^2,\dots,x_k^2)$;
\item
if $\s_\b\in H^{2n}(G_{k}(\C^\infty))$ is expressed as a polynomial $\s_\b=f(c_1,\dots,c_{k})$, then
$$\aligned
\tsaR=&f(p_1,\dots,p_k),\\
\widetilde{\s}_{\a',\R}=&e_{2k}f(p_1,\dots,p_k).
\endaligned$$
\endroster
\endcorollary

\proof A straightforward consequence of
Proposition \ref{Fuks-property}  and Proposition \ref{real-root-map}.
\endproof

\corollary\label{real-schur-vandermond}
Assume that $\a=\b(2)$, $\b=2\g+\bold{r}$, where $\g$ is any $k$-partition and
$\bold{r}=(r,\dots,r)$ is a $k$-partition with $r$ equal to either $0$ or $1$.
Then,
$$s_{\a,\R}(x)=\frac{V_{\b+2\delta}}{V_{2\delta}}=\frac{\vmatrix x_1^{\b_1+2(k-1)}&x_1^{\b_2+2(k-2)}&\dots&x_1^{\b_k}\\
&\dots&\\ x_k^{\b_1+2(k-1)}&x_k^{\b_2+2(k-2)}&\dots&x_k^{\b_k}\endvmatrix}{\vmatrix x_1^{2(k-1)}&x_1^{2(k-2)}&\dots&1\\
&\dots&\\ x_k^{2(k-1)}&x_k^{2(k-2)}&\dots&1\endvmatrix}.
$$
\endcorollary

\proof It follows from Corollary \ref{real-complex-Schur} and
Proposition \ref{schur-vandermond}, since $V_{2\gamma}(x)=V_\gamma(x^2)$
and $V_{2\gamma + \bold{r}}(x)=V_\gamma(x^2) (x_1\dots x_k)^r$ for any $k$-partition $\g$.
\endproof

\lemma\label{basis-orientable}
The classes $\tsaR$ for all even and odd $2k$-partitions $\a$
form an additive basis
of the group $H^*(\tGR{2k}{\infty})/\Tors$.
The corresponding real Schur polynomials $\saR$ form an additive basis in $\Z^{EP}[x]$.
\endlemma

\proof
Classes $\tsaR$ for all even $2k$-partitions $\a$, being dual to a basis in $H_*(\GR{2k}{\infty}$, give an additive basis in $\pii^*(H^*(\GR{2k}{\infty}/\Tors))\subset H^*(\tGR{2k}{\infty})/\Tors)$.
This subgroup is generated by products of the Pontryagin classes $p_i$, by Theorem \ref{EP-ring}. By the same theorem,
$H^*(\tGR{2k}{\infty})/\Tors$ is generated by $1$ and the Euler class $e_{2k}$ as a module over
$H^*(\GR{2k}{\infty}/\Tors)$, so, the classes  $\tsaR$ for all odd $2k$-partitions $\a$,
together with  $\tsaR$ for even $2k$-partitions $\a$ form a basis in $H^*(\tGR{2k}{\infty})/\Tors)$.
\endproof

\subsection{The duality between the classes $[\tCaR]$ and $\tsaR$ in $\tGR{2k}{\infty}$}
In $\tGR{2k}{\infty}$, we have a pair of cells, $\CaR^\pm$, that form the pull-back
of $\CaR$ with respect to the double covering $\pi : \tGR{2k}{\infty}\to \GR{2k}{\infty}$.
The closure of the union $\tCaR=\CaR^+\cup \CaR^-$ yields an integral class of infinite order
$[\tCaR]\in H_{|\a|}(\tGR{2k}{\infty})/\Tors$ if $\a$ is an even or odd $2k$-partition.
Namely, if $\a$ is an even $2k$-partition, we equip $\CaR^\pm$ with the pull-back orientation, so that
$[\tCaR]=\pii^!([\CaR])$. If $\a$ is an odd $2k$-partition, to get an integral class one needs to take the
difference $[\tCaR]=[\CaR^+-\CaR^-]$, where the orientations of $\CaR^\pm$ are chosen to be invariant under the deck transformation.
The latter convention determines the integral class $[\tCaR]$ only up to sign (which is enough for our purpose).

\remark{Remark} These classes $[\tCaR]$ form a basis in $H_*(\tGR{2k}{\infty};\Q)$.
Furthermore, they are divisible by $2$ in the integral homology,
and their halves form an additive basis in $H_*(\tGR{2k}{\infty})/\Tors$.
This can be derived from Theorem 1 in \cite{P},
which treats indeed a more subtle case of finite Grassmannians $\tGR{k}{k+m}$.
\endremark

\lemma\label{Euler-cap}
Consider a Schubert class $[\tCaR]\in H_{|\a|}(\tGR{2k}{\infty})/\Tors$, where $\a=(\a_1,\dots,\a_{2k})$ is an even or odd $2k$-partition. Then
the class $[\tCaR]\cap e_{2k} \in H_{|\a|-2k}(\tGR{2k}{\infty})$
is represented, up to sign and torsion homology elements, by the Schubert class
$\tilde C_{\a',\R}$, where $\a'=\a-\bold{1}=(\a_1-1,\dots,\a_{2k}-1)$, if $\a_{2k}\ge1$.
 If $\a_{2k}=0$, then $[\tCaR]\cap e_{2k}=0$.
\endlemma

\proof
The Euler class $e_{2k}$ restricted to $\tGR{2k}{2k+m}$
is known to be dual to a fundamental class $[\tGR{2k}{2k+m-1}]$ of $\tGR{2k}{2k+m-1}$.
Therefore, the cap product with $e_{2k}$ is realized by the intersection of $[\tCaR]$
with $[\tGR{2k}{2k+m-1}]$  for $m$ sufficiently big with respect to $\vert\a\vert$, which corresponds to subtraction $\bold{1}$ from $\a$.
\endproof

\proposition\label{duality}
Assume that $\a$ and $\b$ are $2k$-partitions, each one is either even or odd, and $|\a|=|\b|=n$.
Then $[\tCaR]\cap\widetilde{\s}_{\b,\R}$ is equal to $\pm2$ if $\a=\b$, and $0$ if $\a\ne\b$.
\endproposition

\proof
For even $\a$, both $[\tCaR]$ and $\tsaR$ are in the invariant subspaces of the covering involution, acting
in the integral homology and cohomology, respectively. For odd $\a$ these classes are both in the skew-invariant
subspaces. This implies that $[\tCaR]\cap\widetilde{\s}_{\b,\R}=0$ if $\a$ and $\b$ have different pairity.

If both $\a$ and $\b$ are even, then $[\tCaR]=\pii^![\CaR]$ and
$\widetilde{\s}_{\b,\R}=\pii^!{\s}_{\b,\R}$. Thus,
$[\tCaR]\cap\widetilde{\s}_{\b,\R}=2([\CaR]\cap{\s}_{\b,\R})=2$
 if $\a=\b$ and $0$ otherwise,
since $\saR$ was defined as the dual to $[\CaR]$.

If both  $\a$ and $\b$ are odd, then $\a=\a'+\bold{1}$, $\b=\b'+\bold{1}$, where $\a'$ and $\b'$ are even.
It follows that
$$[\tCaR]\cap\widetilde{\s}_{\b,\R}=[\tCaR]\cap(e_{2k}\cup\widetilde{\s}_{\b',\R})=
([\tCaR]\cap e_{2k})\cap\widetilde{\s}_{\b',\R}=\pm \tilde C_{\a',\R}\cap\widetilde{\s}_{\b',\R}
$$
where the last identity is due to Lemma \ref{Euler-cap}. So, this case is reduced to the previous one.
\endproof

\subsection{Calculation of the real Schur coefficients}
By Lemma \ref{basis-orientable},
each class $h\in H^*(\tGR{2k}{\infty})/\Tors$
can be decomposed into an integer linear combination
$$h=\sum_{\smallmatrix\text{even and odd}\\ 2k\text{-partitions } \a\endsmallmatrix }\l_\a\tsaR.$$
This  gives the corresponding decomposition of the root polynomial $f=\phi_\R^*h\in \Z^{EP}[z]$:
$$f(x)=\sum_{\smallmatrix\text{even and odd}\\ 2k\text{-partitions } \a\endsmallmatrix } \l_\a\ssaR.$$
The coefficients $\l_\a=\l_a(h)\in\Z$ will be called {\it real Schur coefficients} of $h$,
and of its root polynomial $f$.
If $h\in H^*(\tGR{2k}{\infty})$, we put $\l_\a(h)=\l_\a(h/\Tors)$.

\lemma\label{schur-coefficients}
Assume that $\a$ is an even or odd $2k$-partition, and $h\in H^n(\tGR{2k}{\infty})$.
Then $\l_\a(h)\ovex\frac12[\tCaR]\cap h$ (and in particular, $\l_\a(h)=0$ if $|\a|\ne n$).
\endlemma

\proof By Proposition \ref{duality}, $\{\tsaR\}$ and $\{\frac12[\tCaR]\}$ form
dual, up to sign, basses  in $H_*(\tGR{2k}{\infty};\Q)$.
\endproof

\corollary\label{grassmann-schur-coefficient}
For any class $h\in H^*(\tGR{2k}{\infty})$, we have
$$
\frac12h[\tGR{2k}{m+2k}]\ovex\l_{\bold{m}}(h)\in\Z,
$$
where $\bold{m}$ is the constant $2k$-partition $(m,\dots,m)$.
\endcorollary

\proof
We apply Lemma \ref{schur-coefficients} to
the Grassmannian $\tGR{2k}{m+2k}$, which is a Schubert variety for $2k$-partition $\bold{m}$.
\endproof

Like in Lemma \ref{schur-coefficient}, we can find the coefficients $\l_\a$ using a residue formula.

\lemma\label{real-shur-coefficient}
For any $f\in \Z^{EP}[x]$, $x=(x_1,\dots,x_k)$, and an even or odd $2k$-partition $\a$, we have
$$
\l_\a(f)\ovex\frac1{k!(2\pi i)^k}\int_{T^k}f(x)\overline{\ssaR}(x)V_{2\d}(x)\overline{V_{2\d}}(x)\,\frac{{\roman d}x}{x}
$$
where $T^k=\{x\in\C^k\,: \,|x_1|=\dots=|x_k|=1\}$.
\endlemma

\proof
The proof is the same as that of Lemma \ref{schur-coefficient}, except that we use an expression of the real
Schur polynomials in Corollary \ref{real-schur-vandermond}
instead of the complex one, which was given in Proposition \ref{schur-vandermond}.
\endproof

In the case of $\l_{\bold{m}}$ that we are interested in, we obtain a formula analogous to the one in Corollary \ref{fund-class-evaluation}.

\corollary\label{real-fund-class-evaluation}
For any class $h\in H^{2km}(\widetilde{G}_{2k}(\R^{\infty}))$ its value
on the fundamental class
$[\widetilde{G}_{2k}(\R^{m+2k})]\in H_{2km}(\widetilde{G}_{2k}(\R^{\infty}))$ is determined by the following integral
$$\frac12h\cap[\widetilde{G}_{2k}(\R^{m+2k})]\ovex\l_{\bold{m}}(\phi_\R^*(h))\ovex
\frac1{k!(2\pi i)^k}\int_{T^k}\frac{\phi_\R^*(h)(x)}{x^{\bold{m}}}V_{2\d}(x)\overline{V_{2\d}}(x)\,\frac{{\roman d}x}{x},$$
where $\phi_\R^*(h)\in\Z^{EP}[x]$ is the real root polynomial of $h$.
\qed\endcorollary


\section{Proof of Theorem \ref{main-theorem}}

\subsection{Multivariate integral formula for the number of real $3$-planes on hypersurfaces}
Throughout this section $d\ge1$ is odd, as it is in Theorem \ref{main-theorem}.
We start with a real analogue of the root factorization formula (\ref{root-factorization})
that takes the following form,
in which we call it the {\it real root factorization formula}.

\proposition\label{real-root-formula(square)}
Let $f_d(x)$ denote the real root polynomial of the Euler class
$e_N=e_N(\Sym^d(\til{\tau}^*_{2k,\infty}))\in H^{N}(\tGR{2k}{\infty})$, $N=\binom{d+2k-1}{2k-1}$.
 Then
 $$(-1)^{\frac{N}2}f_d^2(x)=\prod_{\smallmatrix \ell_1+\ell_{\bar 1}+\dots+\ell_k+\ell_{\bar k}=d \\ \ell_i,\ell_{\bar i}\ge0\endsmallmatrix}
 ((\ell_1-\ell_{\bar 1})x_1+\dots+(\ell_k-\ell_{\bar k})x_k).
$$
\endproposition

The factors in the above product formula will be called {\it real root factors}.

\remark{Remark}\label{sign-indeterminacy}
Since we did not fix an orientation of $\Sym^d(\til{\tau}^*_{2k,\infty})$, the polynomial $f_d$
is well-defined only up to sign. In what follows (see formula (\ref{formula-for-f_d}))
we will determine the sign so that the corresponding Euler number,
which we denote $\Cal N_d^e$, becomes positive.
\endremark

\proof
The operation of taking $d$-th symmetric power defines the maps
$\phi_d^\R\:\tGR{2k}{\infty}\to\tGR{N}{\infty}$ and
$\phi_d^\C\:\GC{2k}{\infty}\to\GC{N}{\infty}$
commuting (up to homotopy) with the tautological embedding maps $\o_{2k}$ and $\o_N$
$$
\CD
\tGR{2k}{\infty} @>\phi_d^\R>>\tGR{N}{\infty} \\
@V\o_{2k} VV  @VV\o_N V\\
\GC{2k}{\infty} @>\phi_d^\C>> \GC{N}{\infty}.
\endCD
$$
Then, the result follows from Lemma \ref{real-root-pullback} and Proposition \ref{f_d-formula},
since $e^2(\til\tau^*_{N,\infty})=p_{\frac{N}2}(\til\tau^*_{N,\infty})=
(-1)^{\frac{N}2} \o_N^*(c_{N}(\tau^*_{N,\infty}))$.
\endproof

Theorem \ref{main-theorem} concerns the case $k=2$, and in this case
 Proposition \ref{real-root-formula(square)} reads as follows
 $$f_d^2(x)=\prod_{\smallmatrix \ell_1+\ell_{\bar 1}+\ell_2+\ell_{\bar 2}=d \\ \ell_i,\ell_{\bar i} \ge 0\endsmallmatrix }
 ((\ell_1-\ell_{\bar 1})x_1+(\ell_2-\ell_{\bar 2})x_2).
\eqtag\label{fd-formula}$$

\theorem\label{f_d-counts-N}
Assume that $X\subset P^{m+3}$ is a generic real hypersurface of odd degree $d$, where $\binom{d+3}{3}=4m$.
Then the number of real $3$-planes in
$X$ is finite and bounded from below by
$$
\Cal N_d^e\ovex
\frac1{2(2\pi i)^2}\int_{T^2}\frac{f_d(x)}{x^{\bold{m}}}V_{2\d}(x)\overline{V_{2\d}}(x)\,\frac{{\roman d}x}{x},
$$
where $T^2=\{x\in\C^2\, : \,|x_1|=|x_2|=1\}$
and $f_d(x)$ is the polynomial determined by the formula (\ref{fd-formula}).
\endtheorem

 \proof It follows from Corollary \ref{real-fund-class-evaluation} and Proposition
 \ref{real-root-formula(square)}.\endproof

\subsection{Factorization of $f_d$ in the case of $3$-planes}
The polynomial $g_d(x)=f^2_d(x) \in\Z^{EP}[x]$ can be written as
$$g_d=g_{d,0}g_{d,1}\dots g_{d,\frac{d-1}2},$$ where $g_{d,i}$, $0\le i\le\frac{d-1}2$, is
the product of those real root factors that give
$|\ell_1-\ell_{\bar 1}|+|\ell_2-\ell_{\bar 2}|=d-2i$ (or, equivalently,
$\min(\ell_1,\ell_{\bar 1})+\min(\ell_2,\ell_{\bar 2})=i$).

\lemma We have $g_{d,i}=g_{d-2i,0}^{i+1}$, and thus
$g_d=g_{d,0}g_{d-2,0}^2g_{d-4,0}^3\dots$.
\endlemma

\proof The real root factors in $g_{d,i}$ are the same as in $g_{d-2i,0}$,
but appear as many times as there are partitions
$i=s_1+ s_2$, with $s_j=\min\{\ell_j,\ell_{\bar j}\}\ge0$, $j=1,2$.
\endproof

It is convenient to group the real root factors of $g_d$ by letting
$$
h_d(x_1,x_2)=\prod_{\ell_1,\ell_2\ge1,\ell_1+\ell_2=d}(\ell_1x_1+ \ell_2x_2)
$$
and
$$
\tilde h_d(x_1,x_2)=
\prod_{\ell_1,\ell_2\ge0,\ell_1+\ell_2=d}(\ell_1x_1+ \ell_2x_2)=d^2\,x_1x_2\,h_d(x_1,x_2).
$$

\lemma\label{grouped}
The following identity holds for any odd $d\ge1$:
$$g_{d,0}=
[d^2\,x_1x_2\,h_d(x_1,x_2)h_d(x_1,-x_2)]^2=[\frac1{d^2\,x_1x_2\,}\tilde h_d(x_1,x_2)\tilde h_d(x_1,-x_2)]^2.$$
\endlemma

\proof
The real root factors of $g_{d,0}$ that are taken from $h_d(\pm x_1,\pm x_2)$
go in groups of four elements,
respectively to the four combinations of signs before the coefficients, if the both coefficients do not vanish
(which corresponds in (\ref{fd-formula}) to the cases when both $|\ell_1-\ell_{\bar 1}|$ and $|\ell_2-\ell_{\bar 2}|$ are non-zero and
$\ell_1+\ell_{\bar 1}+\ell_2+\ell_{\bar 2}=d$). Their product converts then to $[h_d(x_1,x_2)h_d(x_1,-x_2)]^2$.
The additional term $d^2\,x_1x_2$ is involved because
if just one of the coefficients $l_1,\bar l_1,l_2,\bar l_2$ is non-zero (and hence, is equal to $d$), then
only two combinations of the signs give root factors contributing to $g_{d,0}$.
\endproof

\corollary
The polynomial $g_d$ is a complete square of the polynomial
$$\align
f_d(x)\ovex&
[d^2\,x_1x_2\,h_d(x_1,x_2)h_d(x_1,-x_2)]
[(d-2)^2\,x_1x_2\,h_{d-2}(x_1,x_2)h_{d-2}(x_1,-x_2)]^2
\dots\\
\ovex&
[d!!(d-2)!!\dots]^2
(x_1x_2)^{\frac{d^2-1}8}\prod_{i=0}^{\frac{d-1}2}[h_{d-2i}(x_1,x_2)h_{d-2i}(x_1,-x_2)]^{i+1}
\qed\endalign
$$
\endcorollary

\example\label{example-degree3}
If $d=1$, then $h_1=1$ (there are no suitable real root factors) and  $g_1=g_{1,0}=x_1^2x_2^2$.
In the first nontrivial case, $d=3$, we have
$$\align h_3(x_1,\pm x_2)=&(2x_1\pm x_2)(x_1\pm 2x_2)=(\pm5x_1x_2+2(x_1^2+x_2^2)),\\
g_{3,0}(x_1,x_2)=&[(9x_1x_2)(5x_1x_2+2(x_1^2+x_2^2))(-5x_1x_2+2(x_1^2+x^2_2))]^2,\\
g_3(x_1,x_2)=&g_{3,0}g_{1,0}^2=[9x^3_1x_2^3(4(x_1^2+x_2^2)^2-25x_1^2x_2^2)]^2,\\
\pm f_3(x_1,x_2)=&9x_1^3x_2^3(4(x_1^2+x_2^2)^2-25x_1^2x_2^2).
\endalign$$
Applying Theorem \ref{f_d-counts-N} and using that
$V_{2\d}\bar V_{2\d}=-\frac{(x_1^2-x_2^2)^2}{x_1^2x_2^2}$,
we conclude that it is the coefficient at $(x_1x_2)^7$ in $-\frac{1}2 f_3(x_1^2-x_2^2)^2$
that gives us the signed count
of real $3$-planes on a generic real $7$-dimensional cubic hypersurface.
This coefficient is equal to $-\frac12(9\times42)=-189$, which gives us $\Cal N_3^e=189$.
 An alternative way to get this coefficient is
to apply Corollary \ref{real-schur-vandermond}
to get the following expressions for the real Schur polynomials, $s_{(5,5,5,5),\R}=(x_1x_2)^5$ and
$$s_{(7,7,3,3),\R}=\frac{V_{(7,3)+(2,0)}}{V_{(2,0)}}=
\frac{\vmatrix x_1^9 & x_1^3\\ x_2^9 & x_2^3\endvmatrix}{\vmatrix x_1^2 & 1\\ x_2^2 & 1\endvmatrix}
=(x_1x_2)^3(x_1^4+x_1^2x_2^2+x_2^4),
$$
which allows us to get the decomposition
$ f_3(x_1,x_2)=9[4s_{(7,7,3,3),\R}-21s_{(5,5,5,5),\R}]$ with
 $\l_{5,5,5,5}=-189$
 and gives us the same result $\Cal N_3^e=189$.

Finally, we can write $\Cal N_3^\R\ge \Cal N_3^{\R,\min}\ge \Cal N_3^e=189$, which implies that a real $7$-dimensional cubic has generically at least
$189$ real $3$-planes.
(See also \cite{FK}, Subsection 6.2, for one more method to find $\Cal N_3^e=189$ by a direct calculation
of the characteristic number $e_{20}[\tGR{4}{9}]$, where $e_{20}=e_4^3(4p_1^2-25e_4^2)$ is represented by the root
polynomial $f_3(x)$).
To compare with, one can find $\Cal N_3^\C=321489$.

Further computations
(using the program Macoley2) show that $\Cal N_5^e= 37655727525$, whereas
$\Cal N_5^\C=64127725294951805931404297113125$.
\endexample

As in Lemma \ref{grouped}, the root factors $(\ell_1x_1+\ell_2x_2)(\ell_2x_1+\ell_1x_2)$ in $h_{d-2i}(x_1,x_2)$
can be grouped together with the factors
$(\ell_1x_1-\ell_2x_2)(\ell_2x_1-\ell_1x_2)$ in $h_{d-2i}(x_1,-x_2)$
to give us the product
$$\align
(\ell_1^2x_1^2-\ell_2^2x_2^2)(\ell_2^2x_1^2-\ell_1^2x_2^2)=&
\ell_1^2\ell_2^2(x_1^2+x_2^2)^2-(\ell_1^2+\ell_2^2)^2x_1^2x_2^2\\
=&(\ell_1\ell_2x_1x_2)^2[(\frac{x_1}{x_2}+\frac{x_2}{x_1})^2-(\frac{\ell_1}{\ell_2}+\frac{\ell_2}{\ell_1})^2]
\eqtag\label{h-factors}
\endalign$$

This product enters in $f_d$ in the power $(i+1)$.
Therefore, we can
put $$f_d(x)=[d!!(d-2)!!\dots]^2(x_1x_2)^{\frac{d^2-1}8}
\prod_{\ell\in\Cal L_d}
[\ell_1^2\ell_2^2x_1^2x_2^2(-\frac{(x_1^2+x_2^2)^2}{x_1^2x_2^2}+\frac{(\ell_1^2+\ell_2^2)^2}{\ell_1^2\ell_2^2})]^{i+1}
\eqtag\label{formula-for-f_d}
$$
where
$$\Cal L_d=\{(\ell_1,\ell_2,i)\,|\,\ell_1+\ell_2=d-2i,\ell_1,\ell_2\ge1,0\le i\le\frac{d-1}2\}.
$$

 Note that in the initial definition of $f_d$ it was defined only up to sign (see
 the remark after Proposition \ref{sign-indeterminacy}),
 and from now on we eliminate this ambiguity by prescribing to $f_d$ the sign given by formula (\ref{formula-for-f_d}).

We can also summarize the above formulae
as follows:
$$\align
h_d(x_1,x_2)h_d(x_1,-x_2)\ovex&\prod_{\smallmatrix \ell_1+\ell_2=d \\ \ell_1,\ell_2\ge1\endsmallmatrix}
(x_1x_2)^2(\ell_1\ell_2)^2[-\frac{(x_1^2+x_2^2)^2}{x_1^2x_2^2}+\frac{(\ell_1^2+\ell_2^2)^2}{\ell_1^2\ell_2^2}], \\
f_d(x) = C_dx_1^Nx_2^N&\prod_{(\ell_1,\ell_2,i)\in\Cal L_d}(\ell_1\ell_2)^{2(i+1)}[-\frac{(x_1^2+x_2^2)^2}{x_1^2x_2^2}+\frac{(\ell_1^2+\ell_2^2)^2}{\ell_1^2\ell_2^2}]^{i+1},
\eqtag\label{h-and-f}\endalign$$
where $C_d=[d!!(d-2)!!\dots]^2$ and $N=\frac12\binom{d+3}3$.

\subsection{Positivity and the maximum}
\proposition\label{max-Fd}
The function $F_d(x)=\frac{f_d(x)}{(x_1x_2)^m}$, where $m=\frac14\binom{d+3}3$ {\rm (}cf. Theorem \ref{f_d-counts-N}{\rm )}
takes positive real values {\rm (}and, in particular, does not vanish{\rm)} at each point of $T^2$.
The maximal value of $F_d(x)$ on $T^2$ is achieved
along the two-component curve $x_2=\pm ix_1$, and this value equals
$$
M_d=C_d\prod_{(\ell_1,\ell_2,i)\in\Cal L_d}(\ell_1^2+\ell_2^2)^{2(i+1)}=
\prod_{i=0}^{\frac{d-1}2}\prod_{\ell=0}^{\frac{d-1}2-i}(\ell^2+(d-2i-\ell)^2)^{2(i+1)}.
\eqtag\label{maximum}
$$
\endproposition

\proof
According to (\ref{h-and-f}), the factors of $(x_1x_2)^{2(1-d)}h_d(x_1,x_2)h_d(x_1,-x_2)$ are equal to
$$
\ell_1^2\ell_2^2[-(t+t^{-1})^2+(k+k^{-1})^2],\ \ \  \text{ where }\ t=\frac{x_1}{x_2},|t|=1,
\text{ and }k=\frac{\ell_1}{\ell_2}, k>0.
$$
 Since $t+t^{-1}$ is real for $|t|=1$, these factors are real, and thus $F_d(x)$ is real for all $x\in T^2$.
Since $0\le(t+t^{-1})^2\le 4\le (k+k^{-1})^2$,
$F_d(x)$ can only vanish for $k+k^{-1}=2$, that is
for $\ell_1=\ell_2$, which is impossible under our assumption that
$d$ is odd.
 The maximal value $(k+k^{-1})^2$ of the factor $\ell_1^2\ell_2^2|-(t+t^{-1})^2+(k+k^{-1})^2|$
is achieved as $(t+t^{-1})^2$ takes its minimal value (equal to $0$), that is along the circles $x_2=\pm ix_1$.

Since these circles are common for all the partitions $(\ell_1,\ell_2)$ involved in the formula for $F_d$,
the product of all the factors achieves its maximum value along the same circles,
and this value is as indicated in the
statement.
\endproof

\proposition\label{asymptotic}
The sequence $M_d$ given  by {\rm (\ref{maximum})} has the following asymptotic grows in the $\log$-scale{\rom :}
$$
\log M_d=\frac1{12}d^3\log d+O(d^3).
$$
\endproposition

\demo{Proof} According to the first order Euler-Maclaurin formula,
the following relations hold for any function $f\in\Cal C^1[0,d]$
and any $r\in\N$:
$$\sum_{\ell=0}^rf(\ell)=\int_0^rf(t)\,{\roman d}t+\frac12[f(0)+f(r)]+\rho_1(r),$$
where $ \rho_1(r)=\int_0^r (x-[x]-\frac12)|f'(t)|\,dt$ and, hence,
$$
|\rho_1(r)|\le\frac1{2}\int_0^r|f'(t)|\,{\roman d}t.
$$
We apply the latter bound
to $f(t)=\log(t^2+(d-t)^2)$ and $r=\frac{d-1}2$.
Since $f'<0$ on $[0,r]$, this gives
$$
|\rho_1(r)|\le\frac1{2}(f(0)-f(\frac{d-1}2))<\log d.
$$

The logarithm of the product
 $J_d=\log(\prod_{\ell=0}^{\frac{d-1}2}(\ell^2+(d-\ell)^2))$
can be written in the form
$$\align
J_d=\sum_{\ell=0}^{\frac{d-1}2}\log(\ell^2+(d-\ell)^2)
=&\int_0^{\frac{d-1}2}\log(t^2+(d-t)^2)\,{\roman d}t\\
+&\log d+\frac12\log(\frac{d^2+1}2)+\rho_1(r)
\endalign$$
where in its turn
$$\aligned
\int_0^{\frac{d-1}2}\log(t^2+(d-t)^2)\,{\roman d}t&=t\log(t^2+(d-t)^2)|_0^{\frac{d-1}2}-\int_0^{\frac{d-1}2} t\frac{2t+2(t-d)}{t^2+(d-t)^2}\,{\roman d}t\\
&=\frac{d-1}2\log(\frac{d^2+1}2)-\tau(d).
\endaligned$$
By $\tau(d)$ we denoted here the integral that can be evaluated as follows:
$$\aligned
\tau(d)=&\int_0^{\frac{d-1}2}t \frac{2t+2(t-d)}{t^2+(d-t)^2}\,{\roman d}t=
\int_0^{\frac{d-1}2}[2+\frac{\frac{d}2(4t-2d)}{2t^2-2dt+d^2}-\frac{d^2}{2t^2-2dt+d^2}]\,{\roman d}t\\
=&(d-1)+\frac{d}2\log(t^2+(d-t)^2)|^{\frac{d-1}2}_0-d\arctan(\frac{2t}d-1)|^{\frac{d-1}2}_0\\
=&(d-1)+\frac{d}2[\log(\frac{d^2+1}2)-2\log d]+d(\arctan{\frac1d}-\frac{\pi}4).
\endaligned$$
This implies a uniform estimate
$\vert \tau(d)\vert \le md+M$
with the constants $m$, $M$ independent of $d$.
Thus, we can write $J_{d}=d\log d+O(d)$, $O(d)\le md+M$, and make an estimate
$$
J_{d-2i}=(d-2i)\log(d-2i)+O(d-2i),\ \ \ O(d-2i)\le md+M.
$$
Then
$$
2(i+1)J_{d-2i}=\log(\prod_{\ell=0}^{\frac{d-1}2-i}(\ell^2+(d-2i-\ell)^2)^{2(i+1)})
= 2(i+1)(d-2i)\log(d-2i)+(i+1)O(d),
$$
and a bound $(i+1)O(d)\le O(d^2)$ gives
$$\align\log M_d=&2\sum_{i=0}^{\frac{d-1}2}[(i+1)(d-2i)\log(d-2i)+O(d^2)]
=2\int_0^{\frac{d-1}2}(t+1)(d-2t)\log(d-2t)\,{\roman d}t\\
+&O(d^3)
=2[-\frac23t^3+\frac12(d-2)t^2+dt]\log(d-2t)|_0^{\frac{d-1}2}+K_d+O(d^3)=K_d+O(d^3)
\endalign$$
where
$$\align
K_d=&4\int_0^{\frac{d-1}2}\frac{\frac23t^3-\frac12(d-2)t^2-dt}{2t-d}\,{\roman d}t=
\frac23\int_0^{\frac{d-1}2}(2t^2-\frac{d-6}2t-\frac{d(d+6)}4)\,{\roman d}t\\
-&\frac23\frac{d^2(d+6)}{4}\int_0^{\frac{d-1}2}\frac{{\roman d}t}{2t-d}
=O(d^3)-\frac1{12}d^2(d+6)\log(d-2t)|_0^{\frac{d-1}2}
\endalign
$$
and finally
$$\log M_d=\frac1{12}d^3\log d+O(d^3).\qed$$
\enddemo

\subsection{Asymptotics}
Consider a sequence of functions
$$F_d(x)=\prod_{(\ell,r)\in {\Cal L}_d}\phi_{\ell,r}(x),\ \  x\in T^2,\ \ \text{ where } d\in\N.$$
Assume that $\phi_{\ell,r}(x)>0$ for all $x\in T^2$, $d$ and $(\ell,r)$,
and that the maximum of $\phi_{\ell,r}$ over $T^2$ is equal to $M_{\ell,r}$ and is achieved on a submanifold
$L\subset T^2$,
which is common for all $d$ and $\ell$.
This is the case in our setting, where $\ell=(\ell_1,\ell_2)$ and
$$\align
\phi_{\ell,r}(x_1,x_2)=&[(\ell_1^2+\ell_2^2)^2-\ell_1^2\ell_2^2\frac{(x_1^2+x_2^2)^2}{x_1^2x_2^2}]^{r+1}\\
=&(\ell_1^2+\ell_2^2)^{2(r+1)}[1-\frac{4\ell_1^2\ell_2^2}{(\ell_1^2+\ell_2^2)^2}\sin^2\phi]^{r+1}
\endalign$$
and $M_{\ell,r}=(\ell_1^2+\ell_2^2)^{r+1}$.
 Namely, if $x_j=e^{i\phi_j}$, $\phi_j\in[0,2\pi]$, $j=1,2$, then
 $$
 \frac{(x_1^2+x_2^2)^2}{x_1^2x_2^2}=(\frac{x_1}{x_2}+\frac{x_1}{x_2})^2=
 (e^{i(\phi_1-\phi_2)}+e^{-i(\phi_1-\phi_2)})^2=4\cos^2(\phi_1-\phi_2)=4\sin^2\phi
 $$
where $\phi=\frac{\pi}2-\phi_1+\phi_2$.
Thus,
$(x_1,x_2)=(e^{i\phi_1},e^{i\phi_2})$
is a point of maximum if and only if $\phi=0$
or $\phi=\pi$, so that in our case $L$ consists of two disjoint circles $\phi_1-\phi_2=\pm\frac{\pi}2$.

We are concerned about the asymptotics, in the logarithmic scale, of the integral
$$I_d=\frac1{8(\pi i)^2}\int_{T^2}F_d(x)W(x)\,\frac{{\roman d}x}x$$
for a certain function $W(x)\ge0$ on $T^2$,
namely, for $W(x)=V_{2\delta}(x)\overline{V}_{2\delta}(x)$,
that is $$\aligned
W(x)&=\vmatrix1&x_1^2\\1&x_2^2\endvmatrix \vmatrix1&x_1^{-2}\\1&x_2^{-2}\endvmatrix
=-\frac{(x_1^2-x_2^2)^2}{x_1^2x_2^2}
=-[e^{2i(\phi_1-\phi_2)}+e^{-2i(\phi_1-\phi_2)}-2]\\
&=2(1-\cos2(\phi_1-\phi_2))
=2(1-\sin 2\phi), \ \ \ \phi=\frac{\pi}2-\phi_1+\phi_2.
\endaligned
$$

\proposition\label{asymptotics}
Under the above choice of functions $F_d$ the following
$\log$-scale asymptotic development holds:
$$\log I_d=\frac1{12}d^3\log d + O(d^3).$$
\endproposition

\proof An upper bound
$$\log I_d\le \frac1{12}d^3\log d + O(d^3)
$$
follows from $\log M_d=\frac1{12}d^3\log d + O(d^3)$ (see Proposition \ref{asymptotic}).

By the localization principle,
the lower bound follows from positivity of $F_d$ and $W$. Indeed, let us consider
the tubular neighborhood $U$ of $L$ defined by $\vert\sin 2 \phi\vert\le\frac12$, that is $|\phi|\le\phi_0=\frac{\pi}{12}$.
Then $1\le W(x)\le3$ for $x\in U$, and
$$
\aligned
\log I_d+\log 8\pi^2&\ge \log\int_{U}F_d(x)W(x)\,\frac{{\roman d}x}x\ge \log[\text{Area}(U)\min_UF_d(x)\min_{U}W(x)]\\
&\ge \log(\frac43 \pi^2)
+ \log(\prod_{(\ell_1,\ell_2,r)\in\Cal L_d}
 (\ell_1^2+\ell_2^2)^{2r}[1-\frac{4\ell_1^2\ell_2^2}{(\ell_1^2+\ell_2^2)^2}
 \sin^2\phi_0]^r )\\
 &= \log(\frac43 \pi^2)+
\log M_d +\Cal R
\endaligned$$
where $$\Cal R=\log\prod_{(\ell_1,\ell_2,r)\in\Cal L_d}[1-\frac{4\ell_1^2\ell_2^2}{(\ell_1^2+\ell_2^2)^2}
 \sin^2\phi_0]^r.$$

Finally, note that
$$|\Cal R|\le c\sum_{(\ell_1,\ell_2,r)\in\Cal L_d} r\frac{4\ell_1^2\ell_2^2}{(\ell_1^2+\ell_2^2)^2}\sin^2\phi_0
\le c d|\Cal L_d|= O(d^3),
$$
where $|\Cal L_d|=\binom{d+2}{2}=O(d^2)$ is the cardinality of $\Cal L_d$ and $c>0$ is a constant independent of $d$.
\endproof

\subsection{Proof of Theorem \ref{main-theorem}}\label{Proof-main}
The asymptotic upper bound for $\Cal N_d^\C$ is established in Proposition \ref{complex-upper-bound}.
The positivity of $\Cal N_d^e$
follows from  Theorem \ref{f_d-counts-N} and the positivity of $F_d$ (see Proposition \ref{max-Fd}).
The asymptotic expression for $\Cal N_d^e$
follows from Theorem \ref{f_d-counts-N} and Proposition \ref{asymptotics}.


\section{Concluding remarks}\label{CR}

\subsection{Counting $3$-planes in the complete intersections} The above approach can be also applied to counting $3$-planes in the complete intersections,
except that the formulas become more cumbersome and the asymptotics is difficult to disclose. For simplicity, let us restrict ourselves to the intersections
of cubic hypersurfaces. Recall that for one cubic hypersurface
(see Example \ref{example-degree3})
the Euler class is given by
$$
e_{20}=e(\Sym^3(\til{\tau}^*_{4,\infty}))\ovex 9e^3(25e^2-4p_1^2)\in H^{20}(\tGR4{\infty}).
$$
An intersection $X$ of $r$ generic real cubic hypersurfaces in $P^{m+3}$
contains a finite number of $3$-planes
if the dimension $20r$ of $e_{20}^r$ is equal to $4m=\dim\tGR4{m+4}$ (see, for example, \cite{DM2}),
that is if $m=5r$.
The signed count of real $3$-planes in $X$ gives a number $\Cal N_{3,\dots,3}^e$,
which is the half of the corresponding count of oriented $3$-planes
that is $e_{20}^r[\tGR4{m+4}]$, so we obtain
$$
2\Cal N_{3,\dots,3}^e\ovex 9^re^{3r}(25e^2-4p_1^2)^r[\widetilde{G}_4(\R^{m+4})].
$$
Due to the multiplicative structure of the Euler-Pontryagin ring, the right-hand side is equal to the result
of substitution of the Catalan numbers $C_k$ instead of $t^k$ in the polynomial $9^r(25-4t)^r$. Standard manipulation
with generating functions and radii of convergence shows that the rate of growth of this sequence
is linear in the logarithmic scale:
$$
\log\Cal N_{3,\dots,3}^e\sim 4r\log 3.
$$

To count the number $\Cal N_{3,\dots,3}^\C$ of complex $3$-planes on a generic intersection of cubic hypersurfaces we can use the multivariate integral Cauchy formula
({\it cf., } Corollary \ref{N_d-via-integral}):
$$
\frac1{4!(2\pi i)^4}\int_{T^4}(\prod_{\smallmatrix\ell_1+\dots+\ell_4=3\\ \ell_1, \dots, \ell_4\ge 0\endsmallmatrix}
(\ell_1z_1+\dots+\ell_4z_4)^r) \frac1{z^{\bold{m}}} V_\d(z)\overline{V_{\d}}(z)\,\frac{{\roman d}z}{z}.
$$
Applying to this integral the saddle point version of the Laplace method, namely, by deforming the integration cycle locally keeping the points of the
maximal absolute value of the product (this value is equal to $3^{20}$) but making the values of the product real at each point of a small neighborhood of the locus of maxima, we obtain
$$
\log\Cal N_{3,\dots,3}^\C\sim 20r\log 3\sim 5\log\Cal N_{3,\dots,3}^e.
$$

\subsection{Another enumerative problem with the $3$-planes}
Note that the Schubert cell
$C_{2,2}\subset\GC{4}{2n+4}$ is formed by the projective $3$-planes in $P^{2n+3}$ which intersect $P^{2n-1}\subset P^{2n+3}$ along a line. Therefore, if we choose a generic set
 $S=\{S_1,\dots,S_{2n}\}$ of projective $(2n-1)$-dimensional subspaces $S_i\subset P^{2n+3}$, then
 the number $\NCp{2n}$
 of $3$-planes $L\subset P^{2n+3}$ such that the intersections $L\cap S_i$, $1\le i\le 2n$, are lines, can be found by evaluation of the power of the Schubert class $\s_{2,2}$:
 $$\NCp{2n}=\s_{2,2}^{2n}[\GC{4}{2n+4}].$$

In the real setting, we define the number $\NRp{2n}$
similarly, by counting the real $3$-planes intersecting a generic
set $S$ of real subspaces $S_i\subset P^{2n+3}$ along lines. This number depends on the choice of a generic set $S$,
and we denote by $\Nminp{2n}$ the minimum of $\NRp{2n}$
with respect to all generic choices of such $S$.
 The number $\s_{2,2,\R}^{2n}[\GR{4}{2n+4}]$ can be interpreted as the signed count of the real $3$-planes, and thus
its absolute value, which we denote by $\Nep{2n}$,
provides an estimate
$$\Nep{2n}\le\Nminp{2n}\le\NRp{2n}\le\NCp{2n}.$$
Here, $\s_{2,2,\R}\in H^4(\GR{4}{2n+4})$ is nothing but the Pontryagin class $p_1$.
 The following result shows that for this enumerative problem
 we get once more a fixed proportion in the logarithmic scale between the
 number of real solutions and the number of complex ones.

\theorem\label{main-p1-count}
Numbers $\NCp{2n}$ and $\Nep{2n}$ have the logarithmic asymptotics
$$\aligned\log\NCp{2n}=&2n\log20+o(n), \\
\log\Nep{2n}=&2n\log2+o(n).
\endaligned$$
In particular,
$$
\log\Nep{2n}\sim\frac1{\log_2 20}\log \NCp{2n}.
$$
\endtheorem

\proposition\label{real-Catalan}
The signed count of the real $3$-planes intersecting each of the given $2n$ generic
$(2n-1)$-planes in $\Rp{2n+3}$ along a line gives the following answer:
 $$
 p_1^{2n}(\tau^*_{4,2n})[G_4(\R^{2n+4})]=\frac1{n+1}\binom{2n}{n}.
 $$
 \endproposition

 \proof
In the Pontryagin ring $H^*(G_4(\R^{2n+4}))/\Tors$
we have ({\it cf.}, Proposition \ref{Fuks-property})
 $$
 p_2^{n}(\tau^*_{4,2n})[G_4(\R^{2n+4})]=
 \sigma_{2,2,2,2,\R}^n[G_4(\R^{2n+4})]= \sigma_{2n,2n,2n,2n,\R}[G_4(\R^{2n+4})]=1.
 $$
Therefore, out statement would follow from
 $p_1^{2n}=\frac1{n+1}\binom{2n}{n}p_2^{n}$.

 On the other hand, there is a similar well known identity $c_1^{2n}=\frac1{n+1}\binom{2n}{n}c_2^{n}$
in  $H^*(G_2(\C^{n+2}))$ which follows easily from the Pieri rule.
 So, it is left to refer to Proposition \ref{Fuks-property}, or just to use the ring isomorphism
 between $H^*(G_4(\R^{2n+4}))/\Tors$ and
$$
H^*(G_2(\C^{n+2}))=
\Z [c_1, c_2, \tilde c_1,\dots, \tilde c_{n}]/\{(1+c_1+c_2)(1+\tilde c_1+\dots +\tilde c_n)=1\}
$$
that puts in correspondence $c_i$ and $p_i$, $i=1,2$.
 \endproof

\corollary
Given a generic set $\{S_1,\dots,S_{2n}\}$ of $(2n-1)$-planes
in $\Rp{2n+3}$, there exist at least $\frac1{n+1}\binom{2n}{n}$ real $3$-planes
that intersect each of the $S_i$'s along a line.
 \endcorollary

 Now, let us look for the asymptotical behavior of the number of complex solutions in the same Schubert problem as that in Proposition \ref{real-Catalan}.

 \proposition\label{complex-Catalan} The following asymptotic approximation holds:
 $$
 \log \sigma_{2,2}^{2n}[G_4(\C^{2n+4})]=2n \log 20 + o(n).
 $$
 \endproposition

 \proof
Applying Corollary \ref{fund-class-evaluation} to $h=\s_{2,2}^{2n}$ we obtain
$$\sigma_{2,2}^{2n}[G_4(\C^{2n+4})]=
\frac1{24(2\pi i)^4}\int_{T^4}\frac{s_{2,2}^{2n}}{x^{2n}}\frac{V^2(x)}{x^3}\,\frac{{\roman d}x}{x}
=\frac1{24(2\pi i)^4}\int_{T^4}f(x)g^n(x)\,\frac{{\roman d}x}{x}
$$
where $f(x)=\frac{V^2(x)}{x^3}=V(x)\overline{V(x)}$ is real and non-negative, as well as
$$g(x)=\frac{s_{2,2}^{2}}{x^2}
=s_{2,2}\overline{s_{2,2}}=|s_{2,2}|^2.$$

 Thus, we deal with a kind of Laplace integrals
and can use the following result, which follows for instance from \cite{MF}, Theorem 2.1 (where the asymptotic
is given under the conditions much weaker than the analyticity that we require below).

 \lemma
 Assume that $f$ and $g$ are real analytic functions taking non-negative values on a compact domain $R\subset\R^d$.
 Let  $M$ be the maximal value of $g$ on $R$.
 Assume in addition that $f$ and $g$ do not vanish identically.
 Then,
 $$
 F(n)= \int_{R} f(x)g^n(x) {\roman d}x
 $$
has the following logarithmic asymptotics
$$\lim_{n\to\infty}\frac{\log F(n)}n=\log M.\quad \qed
$$
\endlemma

Passing to the polar coordinates, $x=\exp(i\theta)$, we obtain
 $$
  \frac1{24(2\pi i )^4} \int_{R} f(\exp(i\theta))g(\exp(i\theta))^n\, {\roman d}\theta,
 $$
where $R=[0,2\pi]^4$, and
 the maximal value of $g(x)=|s_{2,2}(x)|^2$ is $20^2$, since the sum of the coefficients of
 $s_{2,2}=x_1^2x_2^2+ \dots + x_1^2x_3x_4+ \dots + 2x_1x_2x_3x_4$ is $20$.
 \endproof

\demo{Proof of Theorem \ref{main-p1-count}}
It follows from Proposition \ref{complex-Catalan} combined with Stirling formula applied to
Proposition \ref{real-Catalan}.
\qed\enddemo

\subsection{Multivariate integral formula for the number of  real $(2k-1)$-planes on hypersurfaces}
Theorem \ref{f_d-counts-N} and its proof extends easily to the case of counting $2k-1$-planes with any $k\in\N$.
\theorem\label{general-count}
Assume that $X\subset P^{m+2k-1}$ is a generic real hypersurface of odd degree $d$ and that
$\binom{d+2k-1}{2k-1}=2km$.
Then the number, $\Cal N^{\R}_d$, of real $(2k-1)$-subspaces in
$X$ is finite and bounded from below by the
number $\Cal N_d^e\ge 0$ that is given by the multivariate integral formula
$$
\Cal N_d^e\ovex
\frac1{k!(2\pi i)^k}\int_{T^k}\frac{f_d(x)}{x^{\bold{m}}}V_{2\d}(x)\overline{V_{2\d}}(x)\,\frac{{\roman d}x}{x},
$$
where $f_d(x)$ is the polynomial
satisfying the formula of Proposition \ref{real-root-formula(square)}.
\endtheorem

Some heuristic arguments applied to the above integral formula suggest a conjecture that 
$\Cal N_d^e$ is non zero for each 
$k$ and any odd $d$ relatively prime to $k$,
and that for $k$ fixed and such $d$ 
growing to infinity the following asymptotic relation holds
({\it cf.} Conjecture \ref{higher-dim-conj})
$$
\log \Cal N_d^e\sim\frac1{2{(2k-1)!}}d^{2k-1}\log d.
$$



\Refs\widestnumber\key{ABCD}


\ref{BT}
\by R.Bott, L.W Tu
\book Differential forms in algebraic topology.
\bookinfo Graduate Texts
in Mathematics, 82. Springer-Verlag, New York-Berlin,
\yr 1982.
\endref\label{BT}


\ref{Br}
 \by E.H.~Brown, JR
 \paper The cohomology of $BSO_n$ and $SO_n$ with integer coefficients
 \jour Proc. AMS
 \vol 85
 \issue 2
 \yr 1982
 \pages 283--288
\endref\label{Br}


\ref{DM}
 \by O.~Debarre, L.~Manivel
 \paper Sur la vari\'et\'e des espaces lin\'eaires contenus dans une intersection compl\'ete
 \jour Math. Ann
\vol 312
\yr1998
 \pages549--574
\endref\label{DM}

\ref{DM2}
\by  O.~Debarre, L.~Manivel
\paper Sur les intersections compl\`etes r\' eelles
\jour CRAS
 \vol 331 (S\'erie I)
\yr 2000
\pages 887--992
\endref\label{DM2}


\ref{FK}
 \by S.~Finashin, and V.~Kharlamov
 \paper Abundance of real lines on real projective hypersurfaces
 \jour Int. Math. Res. Notices
 \yr 2013
 \pages  3639--3646
\endref\label{FK}

\ref{F}
 \by D.~B.~Fuks
 \paper Classical manifolds. (Russian)
 \inbook Current problems in mathematics. Fundamental directions
\vol 12
\yr 1986
\pages 253--314
\endref\label{Fuks}

\ref{GZ}
\by P.~Georgieva and A.~Zinger
\jour arXiv:1309.4079
\paper Enumeration of real curves in $\C\Bbb P^{2n-1}$ and a WDVV relation for real Gromov-Witten invariants
\endref\label{GZ}

\ref{GM}
 \by D.B.~Gr\"unberg, P.~Moree
 \paper Sequences of enumerative geometry: congruences and asymptotics (with Appendix by Don Zagier)
 \jour Experiment. Math.
 \vol 17
 \yr 2008
\pages 409--426
\endref\label{GM}


\ref{IKS}
\by I.~Itenberg, V.~Kharlamov, and E.~Shustin
\paper Welschinger invariant and enumeration of real rational curves.
\jour Int. Math. Res. Not.
\yr 2003
\issue 49
\pages 2639--2653
\endref\label{iks}

\ref{IKS2}
\by I.~Itenberg, V.~Kharlamov, and E.~Shustin
\paper Welschinger invariants of real del Pezzo surfaces of degree $\ge 3$
\jour Math. Ann.
\vol 355
\yr 2013
\issue 3
\pages 849--878
\endref\label{iks2}

\ref{KR}
\by V.~Kharlamov and R.~Rasdeaconu
\paper Counting Real Rational Curves on K3 Surfaces
\jour Int. Math. Res. Notices
\yr 2014
\vol doi: 10.1093/imrn/rnu091
\endref\label{KR}

\ref{MF}
 \by V.P. Maslov and M.V. Fedoryuk
 \paper Logarithmic asymptotic of rapidly decreasing solutions of Petrovskii hyperbolic equations
 \jour Math. Notes
 \vol 45:5
 \yr 1989
 \pages 382--391
\endref\label{MF}

\ref{OT}
\by    C.~Okonek,  A.~Teleman
\paper Intrinsic signs and lower bounds in real algebraic geometry
\jour J. Reine Angew. Math.
\vol 688
\yr 2014
\pages 219--241
\endref\label{OT}

\ref{P}
 \by L.~Pontryagin
 \paper Characteristic cycles on differentiable manifolds. (Russian)
\jour Mat. Sbornik N. S.
\vol 21(63)
\yr 1947
\pages 233--284
\endref\label{P}


\ref{St}
\by R.~P.~Stanley
 \book Enumerative Combinatorics. Volume 2
 \bookinfo Cambridge Studies in Applied Mathematics 62,
 Cambridge University Press, Cambridge
 \yr 1999
\endref\label{St}

\endRefs
\enddocument